\documentclass[reqno,11pt]{article}
\usepackage{amsmath,amssymb}
\usepackage{geometry}
\usepackage{amsmath, amsfonts,amsthm,amssymb,amsbsy,upref,color,graphicx,amscd,enumerate}
\usepackage[active]{srcltx}
\usepackage[utf8]{inputenc}
\usepackage{algorithm}
\usepackage{stmaryrd}
\usepackage{array}
\usepackage{algorithmicx}
\usepackage{algpseudocode}
\usepackage{graphicx}
\usepackage{color}
\usepackage[table,xcdraw,svgnames]{xcolor}
\definecolor{light-salmon}{RGB}{255,140,120}
\graphicspath{{./metapost/}{./pics/}}

\usepackage[colorlinks]{hyperref}
\hypersetup{
	allbordercolors=.,
	citecolor=DarkGreen,
	linkcolor=DarkRed,
	urlcolor=NavyBlue
}

\usepackage{enumitem}

\theoremstyle{plain}
\newtheorem{thm}{Theorem}

\newtheorem{prop}[thm]{Proposition}

\theoremstyle{definition}

\renewcommand{\Bbb}{\mathbb}
\newcommand{\bo}[1]{{\bf #1}}

\makeatletter
\DeclareFontFamily{U}{tipa}{}
\DeclareFontShape{U}{tipa}{m}{n}{<->tipa10}{}
\newcommand{\arc@char}{{\usefont{U}{tipa}{m}{n}\symbol{62}}}%

\newcommand{\arc}[1]{\mathpalette\arc@arc{#1}}

\newcommand{\arc@arc}[2]{%
	\sbox0{$\m@th#1#2$}%
	\vbox{
		\hbox{\resizebox{\wd0}{\height}{\arc@char}}
		\nointerlineskip
		\box0
	}%
}
\makeatother

\title{A Geometric proof for the Polygonal Isoperimetric Inequality}

\author{Beniamin Bogosel}

\begin{document}
	 \maketitle
	 
	 \begin{abstract}
	 	Gradients of the perimeter and area of a polygon have straightforward geometric interpretations. The use of optimality conditions for constrained problems and basic ideas in triangle geometry show that polygons with prescribed area minimizing the perimeter must be regular.
	 \end{abstract}
	 
	 Given a closed curve $\gamma$ in the plane bounding a region $\omega$ having a given area $A_0$, what is the shortest possible value for the length $L$ of $\gamma$? This problem is called the \emph{isoperimetric inequality} and has been widely studied since the antiquity. The solution to the isoperimetric problem is, without surprise, the circle. Historical details and various proof ideas are shown in \cite{Isop-AMM}. For a positive integer $n \geq 3$, the same question can be formulated, restricting to the family of $n$-gons:
	 
	 Which $n$-gon $P$ of area $A(P)$ has the smallest perimeter $L(P)$? 
	 
	 Denoting by $\mathcal P_n$ the class of simple $n$-gons having $n$ vertices, the polygonal isoperimetric inequality can be written in the following form
	 \begin{equation}
	 \min_{\Omega \in \mathcal P_n, A(P)=A_0} L(P).
	 \label{eq:isop-poly}
	 \end{equation}
	 Drawing a parallel between the continuous and the discrete case, it is natural to expect that the solution to problem  \eqref{eq:isop-poly} is the regular $n$-gon. Various proofs showing the optimality of the regular $n$-gon are known, some of which involve area preserving transformations that decrease the perimeter.
	 
	 A proof based on optimality conditions is recalled in \cite{Isop-AMM}. Since any $n$-gon $P \in \mathcal P_n$ depends on $2n$ real variables, namely the coordinates of the vertices, \eqref{eq:isop-poly} is, in fact, a constrained minimization problem of the form 
	 \begin{equation}
	 \min_{ g(x)=0} f(x),
	 \label{eq:min-constr}
	 \end{equation}
	 where $f,g : \Bbb{R}^{2n} \to \Bbb{R}$. Classical theory of optimality conditions for this kind of problems states the following:
	 \begin{multline}
	 \text{If } x^* \text{ solves } \eqref{eq:min-constr} \text{ and }\nabla g(x^*)\neq 0,
	 \text{ then there exists } \lambda \in \Bbb{R} \text{ such that } \nabla f(x^*) = \lambda \nabla g(x^*).
	 \label{eq:optim-constr}
	 \end{multline}
	 This type of results is standard in the theory of optimization and more details can be found in \cite[Section 5.5.3]{cvxoptim} for example. Writing \eqref{eq:optim-constr} using the expressions for the area and perimeter in terms of the vertex coordinates shows that the regular $n$-gon is the only possible critical point. This proof recalled in \cite{Isop-AMM} uses algebraic aspects regarding complex numbers.
	 
	 The purpose of this note is to show the geometric interpretation of \eqref{eq:optim-constr} when applied to solution of the polygonal isoperimetric problem. A geometric description of the gradients of the area and perimeter of an $n$-gon is provided. Then, the optimality condition \eqref{eq:optim-constr} combined with well known facts from triangle geometry give a straightforward argument showing that the regular $n$-gon is the only solution to problem \eqref{eq:isop-poly}.
	 
	 \section{Geometric description of the gradients}
	 
	 Denote by $\bo a_0,...,\bo a_{n-1}$ the vertices of a polygon $P$ oriented in a counter-clockwise direction. For simplicity we assume that the indices are taken modulo $n$. The perimeter and the area of an $n$-gon vary smoothly, as long as two vertices do not merge. This can readily be seen by writing the analytical expressions of the area and perimeter in terms of the vertex coordinates $\bo a_i(x_i,y_i)$, $i=0,...,n-1$:
	 \[ A(P) = \frac{1}{2} \sum_{i=0}^{n-1} (x_iy_{i+1}-x_{i+1}y_i).\]
	 \[ L(P) = \sum_{i=0}^{n-1} \sqrt{(x_{i+1}-x_i)^2+(y_{i+1}-y_i)^2}.\]
	 Therefore, like in \cite{Isop-AMM}, we may compute the gradients, i.e. the vectors containing all partial derivatives for the perimeter and area functionals. However, in the following we focus on giving a geometric interpretation for these gradients. Recall the following properties verified by the gradient of a function. We omit the proofs, which are found in any generic multivariable calculus course like \cite[Section 2.6]{calculus}.
	 
	 \begin{prop}
	  Let $f: \Bbb{R}^d \to \Bbb{R}$ be a differentiable function. The gradient $\nabla f(x)=(\frac{\partial f}{\partial x_i}(x))_{i=1}^d$ at the point $x$ is a vector in $\Bbb{R}^d$ which verifies the following: 
	 
	 (i) $\nabla f(x)$ points in the direction of the \emph{steepest ascent}, i.e. the direction along which $f$ increases the fastest around the point $x$.
	 
	 (ii) The norm of the gradient $|\nabla f(x)|$ is equal to the largest value of the directional derivative $h'(t)$ where $h(t) = f(x+t d)$ and $d$ is an arbitrary unit vector. 
	 
	 (iii) The gradient $\nabla f(x)$ is orthogonal to the level set $\{ y \in \Bbb{R}^d :  f(y) = f(x)\}$ passing through $x$.
	 \label{prop:grad}
	\end{prop}

	 We start by observing that when moving the vertex $\bo a_i$ the variation of the perimeter and area are the same for the polygon $P$ and for the triangle $\Delta \bo a_{i-1}\bo a_i \bo a_{i+1}$. Therefore, to compute the gradients of the perimeter and area we with respect to coordinates of $\bo a_i$ we can work on $\Delta \bo a_{i-1}\bo a_i \bo a_{i+1}$. 
	 
	 \bo{Gradient of the Area.} It is well known that the area of a triangle with fixed base is proportional to the height from the third vertex. As a direct consequence of Proposition \ref{prop:grad} we find that the gradient of the area with respect to coordinates of vertex $\bo a_i$ is a vector $\vec v_i$ such that:
	 
	 (a) $\vec v_i$ is orthogonal to the segment $\bo a_{i-1}\bo a_{i+1}$ and points towards the exterior of $P$.
	 
	 (b) $\vec v_i$ has length equal to $\frac{1}{2} | \bo a_{i-1}\bo a_{i+1}|$, a direct consequence of the area formula for a triangle.
	 
\begin{figure}
	 	\vspace{0.25cm}
	 	
	 	\includegraphics[width=0.45\textwidth]{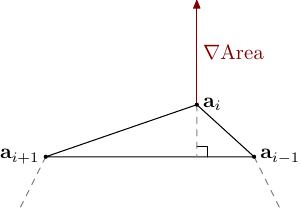}\quad
	 	\includegraphics[width=0.45\textwidth]{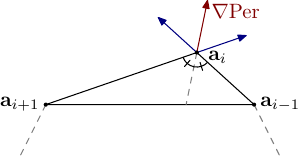}
	 	\caption{Geometric interpretation for the gradient of the area and perimeter of a polygon.}
	 	\label{fig:geometric}
	 \end{figure}

	 Using the notation $\mathcal R_\theta$ for a rotation of angle $\theta$ in the trigonometric direction, we have the more precise description given by:
	 \begin{equation}
	 \vec v_i = \frac{1}{2}\mathcal R_{-\pi/2} (\overrightarrow{\bo a_{i-1}\bo a_{i+1}}).
	 \label{eq:grad-area}
	 \end{equation}
	 An illustration is given in Figure \ref{fig:geometric}. 
	 
	 \bo{Gradient of the Perimeter.} We start by recalling that if $\bo a$ and $\bo o$ are two distinct points in the plane, the gradient of the distance $|\bo o \bo a|$ with respect to $\bo a$ is given by the unit vector pointing from $\bo o$ to $\bo a$ (a direct consequence of Proposition \ref{prop:grad}). Therefore
	 \[ \nabla_{\bo a_i} |\bo a_{i-1}\bo a_i| = \frac{1}{|\bo a_{i-1}\bo a_i|} \overrightarrow{\bo a_{i-1}\bo a_i} \text{ and } \nabla_{\bo a_i} |\bo a_{i+1}\bo a_i| = \frac{1}{|\bo a_{i+1}\bo a_i|} \overrightarrow{\bo a_{i+1}\bo a_i}.\]
	
	 Thus, the gradient of the perimeter of $P$ with respect to $\bo a_i$ is 
	 \begin{equation} \vec w_i = \frac{1}{|\bo a_{i-1}\bo a_i|} \overrightarrow{\bo a_{i-1}\bo a_i}+ \frac{1}{|\bo a_{i+1}\bo a_i|} \overrightarrow{\bo a_{i+1}\bo a_i}.
	 \label{eq:grad-perim}
	 \end{equation}
	 Since the sum of two vectors having the same magnitude is aligned with the bisector of the angle made by the two vectors we have the following:
	 
	 (a) $\vec w_i$ is aligned with the bisector of the angle $\angle \bo a_{i-1}\bo a_i \bo a_{i+1}$. 
	 
	 (b) $\vec w_i$ points towards the exterior when $\theta_i=\angle \bo a_{i-1}\bo a_i \bo a_{i+1}<\pi$ and towards the interior if $\theta_i>\pi$. Moreover, $\vec w_i=\vec 0$ if and only if $\bo a_{i-1},\bo a_i, \bo a_{i+1}$ are aligned in this order.
	 
	 (c) $\vec w_i$ has length equal to $2\cos \frac{\theta_i}{2}$.
	 
	 An illustrative example is shown in Figure \ref{fig:geometric}. 
	 
	 To further illustrate the geometry of the gradients, in Figure \ref{fig:geometric} three situations are plotted, an arbitrary $n$-gon, a star-shaped $n$-gon and the regular one. For each vertex, the corresponding components of the gradient are plotted. Bold arrows represent the gradient of the perimeter. It can be observed that in the regular case, gradients are colinear with the same factor of proportionality for each vertex. In the next section we show that this only occurs when the $n$-gon is regular. For a star-shaped polygon, the gradients have the same direction at a fixed vertex, but the vectors do not scale with the same factor. For arbitrary $n$-gons, the gradients may not even have the same directions for a given vertex.
 
 \begin{figure}
 	\vspace{0.25cm}
 	
 	\centering
 	\includegraphics[height=0.31\textwidth]{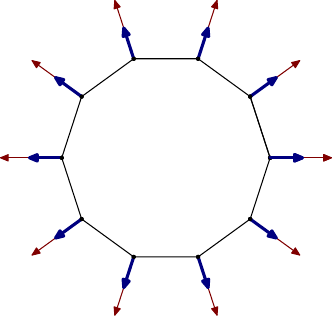}\quad
 	\includegraphics[height=0.31\textwidth]{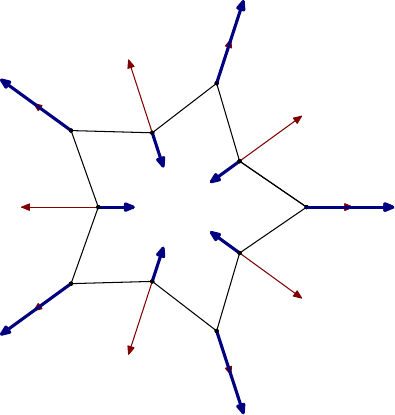}\quad
 	\includegraphics[height=0.31\textwidth]{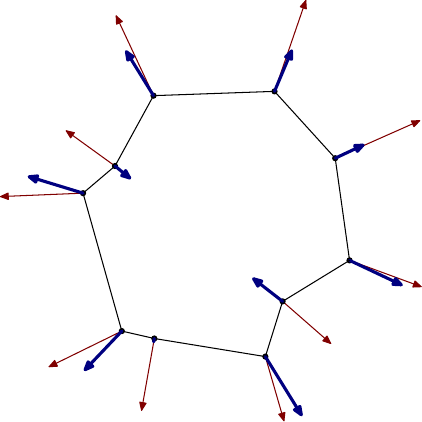}
 	\caption{Gradients of the perimeter and area for regular, star-shaped and arbitrary $n$-gons. Bold arrows represent the gradient of the perimeter.}
 	\label{fig:grad-examples}
 \end{figure}
	 
	 \section{Any optimal polygon must be regular}
	 
	 We are now ready to characterize solutions of \eqref{eq:isop-poly} via optimality conditions. However, before doing this we must make sure that optimal polygons solving \eqref{eq:isop-poly} exist and have some elementary qualitative properties. A proof is given in \cite{Isop-AMM}, but for the sake of completeness, the main ideas are recalled below:
	 \begin{itemize}
	 	\item Taking the convex envelope of a polygon, i.e. replacing the polygon with the smallest convex polygon containing it, decreases its perimeter and increases its area. Therefore, it is not restrictive to consider \eqref{eq:isop-poly} in the class of convex $n$-gons. All arguments below can be tweaked to handle general $n$-gons, but certain aspects are simpler to explain in the convex case.
	 	\item Since $L(P) \geq 0$ there exists a minimizing sequence $P_n$ such that $A(P_n)=A_0$ and $L(P_n)$ converges to $\inf\limits_{P\in \mathcal P_n, A(P)=A_0}L(P)$.
	 	\item The perimeter of $P_n$ is larger than the diameter of $P_n$, therefore the diameters of $P_n$ are bounded. Without restriction of generality, we may suppose that $P_n$ are contained in a closed ball.
	 	\item Classical compactness arguments in $\Bbb{R}^d$, i.e. bounded sequences contain converging sub-sequences, show that there exists a polygon $P^*$ such that, up to a subsequence, $P_n \to P^*$. The area and perimeter are continuous in terms of vertex coordinates, therefore $P^*$ is a solution of \eqref{eq:isop-poly}, provided it has $n$ distinct vertices.
	 	\item $P^*$ is convex and has at most $n$ vertices, but notice that vertices may merge in the limiting process. If $P^*$ has fewer than $n$ vertices we add, repeatedly if necessary, a midpoint of a current edge as a new variable vertex. We may assume, therefore that $P^*$ has $n$ vertices. 
	 \end{itemize}
	 
	Consider $P^*$ an optimal $n$-gon for problem \eqref{eq:isop-poly}. In view of the discussion in the previous section, the gradient of the area is non-zero as soon as there exists $i$ such that $|\bo a_{i-1}\bo a_{i+1}|\neq 0$. This is obviously the case for a minimizer $P^*$, since $|\bo a_{i-1}\bo a_{i+1}|=0$ implies that the area of $\Delta \bo a_{i-1}\bo a_i \bo a_{i+1}$ is zero. The optimality of $P^*$ implies that $\bo a_{i-1},\bo a_i$ and $\bo a_{i+1}$ coincide, which is impossible. Moreover, the previous discussion implies that $\vec v_i\neq 0$, for every $i=0,...,n-1$, where $v_i$ are given by \eqref{eq:grad-area}.
	
	Thus, we conclude that $\nabla A(P^*)$ is not identically zero and the characterization of the optimality conditions given in \eqref{eq:optim-constr} applies. Therefore, there exists $\lambda \in \Bbb{R}$ such that
	\[ \nabla L(P^*) = \lambda \nabla A(P^*).\]
	In other words, for $v_i$ and $w_i$ defined by \eqref{eq:grad-area} and \eqref{eq:grad-perim}, respectively we have $\vec w_i = \lambda \vec v_i$. 
	
	Let us first observe that $\lambda \neq 0$. Indeed, there exist at least three consecutive vertices which are not collinear, otherwise all points would lie on the same line and the area of $P^*$ would be zero. If $\bo a_{i-1}, \bo a_i, \bo a_{i+1}$ are not colinear then the gradient of the perimeter with respect to $\bo a_i$ given by $w_i$ in \eqref{eq:grad-perim} is not zero. Therefore $\nabla L(P^*)\neq 0$, which implies $\lambda \neq 0$. Finally, since $\vec w_i = \lambda \vec v_i$ and $\vec v_i \neq 0$, the gradients of the perimeter $w_i$ given in \eqref{eq:grad-perim} also verify $w_i\neq 0$ for every $i=0,...,n-1$.
	
	Let us now use the three implications of a vectorial equality. 
	
	(a) \bo{Two equal vectors have the same direction.} For each $i=0,...,n-1$ vectors $\vec w_i$ and $\vec v_i$ are non-zero and collinear. Since the first one is aligned with the bisector of $\angle \bo a_{i-1}\bo a_i \bo a_{i+1}$ and the second is aligned with the height from $\bo a_i$ in $\Delta \bo a_{i-1}\bo a_i \bo a_{i+1}$, it follows that this triangle is isosceles with respect to vertex $\bo a_i$. Therefore for every $i=0,...,n-1$ we have $|\bo a_{i-1}\bo a_i| = |\bo a_i\bo a_{i-1}|$, which means that all edges of $P^*$ have the same length $\ell$.
	
	(b) \bo{Two equal vectors have the same orientation.} Since $v_i$ points always towards the exterior of the polygon it follows that $w_i$ will always point towards the exterior. Therefore $P^*$ does not have angles greater than $\pi$. This fact was already established, during the discussion regarding existence of solutions, since optimal polygons are convex.
	
	(c) \bo{Two equal vectors have the same length.} In the isosceles triangle $\Delta \bo a_{i-1}\bo a_i \bo a_{i+1}$ we have $|\bo a_{i-1}\bo a_{i+1}| = 2\ell \sin \frac{\theta_i}{2}$ where $\theta_i = \angle \bo a_{i-1}\bo a_i \bo a_{i+1}$ and $\ell$ is the length of the edges of the polygon. Recalling that the length of $\vec w_i$ is $2\cos \frac{\theta_i}{2}$ we find that
	\[ \frac{1}{2} \ell \sin \frac{\theta_i}{2} = 2\lambda \cos \frac{\theta_i}{2},\]
	for every $i=0,...,n-1$. Therefore $\tan \frac{\theta_i}{2} = \frac{4\lambda}{\ell}$ showing that all angles of $P^*$ are equal.
	
	In conclusion $P^*$ has equal edge lengths and equal angles, therefore $P^*$ is the regular $n$-gon. The Lagrange multiplier $\lambda$ can be explicited in terms of the edge length $\ell$ and the angle $\theta$ associated to the regular polygon of area $A$ using the relations
	\[ \lambda = \frac{1}{4}\ell \tan \frac{\theta}{2},\ A = \frac{1}{4}n\ell^2 \tan \frac{\theta}{2}.\]

	Viewing gradients geometrically, rather than algebraically, turns the optimality condition \eqref{eq:optim-constr} verified by an optimal polygon in \eqref{eq:isop-poly} into a series of elementary observations which show that the optimal $n$-gon must be regular.
	 
\bibliography{./DiscreteIsop}
\bibliographystyle{abbrv}

		\noindent Beniamin \textsc{Bogosel}, Centre de Math\'ematiques Appliqu\'ees, CNRS, École polytechnique, Institut Polytechnique de Paris, 91120 Palaiseau, France\\
		\nolinkurl{beniamin.bogosel@polytechnique.edu}

\end{document}